\documentclass[12pt, reqno]{amsart}
\usepackage{color, MnSymbol, bm, amsmath}
\usepackage{multirow}

\usepackage[authoryear]{natbib}
\usepackage{hyperref}
\usepackage[latin1]{inputenc}
\usepackage{graphicx}
\usepackage{bbm}
\usepackage{subcaption}
\usepackage{epstopdf}

\numberwithin{equation}{section}

\theoremstyle{plain}                    
\newtheorem{thm}{Theorem}[section]
\newtheorem{lem}{Lemma}[section]

\theoremstyle{definition}

\newtheorem{exmp}{Example}[section]

\theoremstyle{remark}

\newtheorem{ack}{Acknowledgment}

\def\N{{\mathbb N}}
\def\R{{\mathbb R}}

\def\Z{{\mathbb Z}}

\def\1{{\mathbbm{1}}}

\def\argmin{\mbox{arg}\,\min}

\newcommand{\cov}{\mathop{\rm cov}\nolimits}

\newcommand{\be}{\begin{equation}}
\newcommand{\bd}{\begin{displaymath}}
\newcommand{\ed}{\end{displaymath}}
\newcommand{\bea}{\begin{eqnarray}}
\newcommand{\eea}{\end{eqnarray}}
\newcommand{\bean}{\begin{eqnarray*}}
\newcommand{\eean}{\end{eqnarray*}}

\begin{document}

\thispagestyle{empty}

\baselineskip10pt

\begin{center}
{\large \bf On Integrated $L^{1}$ Convergence Rate of an Isotonic Regression Estimator
for Multivariate Observations}
\end{center}

\vspace*{0.5cm}

\begin{center}
Konstantinos Fokianos\\
University of Cyprus\\
Department of Mathematics and Statistics\\
P.O.~Box 20537\\
CY -- 1678 Nicosia\\
Cyprus\\
E-mail: fokianos@ucy.ac.cy\\[1cm]
Anne Leucht\\
Technische Universität Braunschweig\\
Institut f\"ur Mathematische Stochastik\\
Universit\"atsplatz 2\\
D -- 38106 Braunschweig\\
Germany\\
E-mail: a.leucht@tu-bs.de\\[1cm]
Michael H.~Neumann\\
Friedrich-Schiller-Universit\"at Jena\\
Institut f\"ur Mathematik\\
Ernst-Abbe-Platz 2\\
D -- 07743 Jena\\
Germany\\
E-mail: michael.neumann@uni-jena.de\\[1.5cm]
\end{center}

\vspace*{10pt}

\begin{center}
{\bf Abstract}
\end{center}
\begin{small}
We consider a general monotone regression estimation where we allow for independent and dependent
regressors. We propose a modification of the classical isotonic least squares estimator
and establish its rate of convergence for the integrated $L_1$-loss function. The methodology
captures the shape of the data without assuming additivity or a parametric form for the regression function.
Furthermore, the degree of smoothing is chosen automatically
and no auxiliary tuning is required for the theoretical analysis.
Some simulations and two real data illustrations complement the study of the proposed estimator.
\end{small}

\vspace*{1.0cm}

{\bf Keywords:} Isotonic least squares estimation, multivariate isotonic regression,
nonparametric estimation, rate of convergence, shape constraints, strong mixing.\\

%
%
%

\newpage

\setcounter{page}{1}
\pagestyle{headings}
\normalsize
\section{Introduction}
\label{S1}

We consider the classical mean regression model
\begin{equation}
\label{eq.motivation}
Y_t=f(I_t)+\varepsilon_t\qquad\text{with}\quad E(\varepsilon_t\mid I_t)=0\;\; a.s., \quad t\in \Z,
\end{equation}
where we assume that the regression function $f\colon\,D \to \R$, $D\subseteq \R^d$, is unknown and allow for both independent and dependent observations
$((Y_t, I_t')')_{t}$ (here and in the sequel, $x'$ denotes the transpose of a vector~$x$).

Notably,  the problem of estimating a regression function subject to shape constraints,
in the context of time series, has not been addressed adequately in the literature, to the best of our knowledge.
There exist  a large body of literature  on estimation
and testing for situations where the class of admissible functions~$f$  can be
parametrized by a finite-dimensional parameter; see e.g.~\citet{Esc06}, \citet{FZ10}  and \citet{SS11} among others.
There are also many results on nonparametric kernel estimators for $f$ relying on the assumption that
the covariate vector~$I_t$ has a Lebesgue density. For an overview, we refer the reader to the
monographs by \citet{H90} and \citet{FG96}. On the other hand, there are
numerous applications that the covariates do not possess a density with respect to the
Lebesgue measure; a case in point is various count time series models which have been
employed for the analysis of financial data (e.g. modeling the number of transactions) or
biomedical data (e.g. modeling infectious diseases); see Fokianos et al. (2009) for instance and Sec. \ref{Sub-Sec:data}.

The primary aim of this  work is to provide integrated $L^{1}$-loss  convergence rate of  a nonparametric
estimator of~$f$ \emph{subject to shape constraints without assuming additivity};
in particular we  assume throughout this work  that the function~$f$
in \eqref{eq.motivation} is isotonic. The assumption of isotonicity seems to be
appropriate in the context of many applications and, in fact, some popular
parametric models share this property,
for example, autoregressive and GARCH type models with nonnegative coefficients.
Application of standard nonparametric methods such as kernel estimators
of the function~$f$, as proposed e.g.~by  \citet{MS94}, \citet{DNP06}, \citet{CFG09}, \citet{DP13} or generalized additive modeling \citet{CS16} among other references,
depends on  a data-driven choice of smoothing parameters, such as a bandwidth.
While the simple leave-one-out cross-validation may fail, the method of
leave-$k$-out cross-validation involves a choice of~$k$, which in turn requires a difficult
subjective decision.

Another popular shape-constrained estimator of the function~$f$
is the isotonic least squares estimator (LSE)~$\widetilde{f}_n$ which is given by
\begin{displaymath}
\widetilde{f}_n \,\in\, \argmin_{g\; \mbox{{\scriptsize isotonic}}} \sum_{t=1}^n \left( Y_t - g(I_t) \right)^2.
\end{displaymath}
In sharp contrast to usual kernel estimators, the isotonic least squares estimator does not require the choice of
any smoothing parameter since an appropriate tuning of the degree of smoothing is done
automatically. This estimator seems to be less sensitive to irregularities in the design
and if the target function is indeed isotonic then this estimator is consistent;
see e.g.~\citet{CT87} and references therein.

Denote by $\1(\cdot)$ the indicator function. Then, it is well known that~$\widetilde{f}_n$ satisfies at all observation points
$x\in\{I_1,\ldots,I_n\}$ the following equations:
\begin{subequations}
\begin{eqnarray}
\label{2.1}
\widetilde{f}_n(x) & = & \max_{U\colon\; x\in U} \; \min_{L\colon\; x\in L} \mbox{Av}_Y(L\cap U) \\
\label{2.2}
& = & \min_{L\colon\; x\in L} \; \max_{U\colon\; x\in U} \mbox{Av}_Y(L\cap U),
\end{eqnarray}
\end{subequations}
where
\begin{displaymath}
\mbox{Av}_Y(B) \,=\, \frac{ \sum_{t=1}^n Y_t \; \1( I_t\in B ) }{ \#\{t\leq n\colon \;\; I_t\in B\} },\qquad B\subseteq [0,1]^d,
\end{displaymath}
and $U$ and $L$ denote upper and lower sets, respectively;
see e.g.~Theorem~1 in \citet{Bru55} and Theorem~1.4.4 in \citet[p.~23]{RWD88}.
(A set $U\subseteq [0,1]^d$ is called an upper set if $x\in U$ and $x\preceq y$ imply that $y\in U$.
Analogously, $L\subseteq \R^d$ is called a lower set if $x\in L$ and $x\succeq y$ imply that $y\in L$. Here,
the notation  $x\preceq y$ ($x\succeq y$, respectively) denotes  that $x_i\leq y_i$  ($x_i\geq y_i$, respectively), for all $i=1,\ldots,d$.)
While~$\widetilde{f}_n$ is uniquely defined at the observation points, there is some arbitrariness  of
choosing~$\widetilde{f}_n$ between these points; only the postulated isotonicity has to be satisfied.

For the univariate case, i.e.~$d=1$, there are already several  results reported in the literature concerning
the asymptotic behavior (usually assuming a deterministic regressor)
of the classical isotonic least squares estimator~$\widetilde{f}_n$. Pointwise asymptotic distributions
of isotonic least squares estimators assuming  short and long
range dependence of the error sequence $(\varepsilon_t)_{t}$ have been derived recently  by \citet{AH06} and \citet{DMP11}.
In particular, it is known that this estimator converges at the optimal rate~$n^{-1/3}$ to~$f$.
\citet[Theorem~2.3]{Zha02} studies the  case of independent but not necessarily identically distributed  errors
and shows that
$(n^{-1}\sum_{i=1}^n E(\widetilde{f}_n(t_i)-f(t_i))^p)^{1/p}=O(n^{-1/3})$,
where $t_1,\ldots,t_n$ are  values of a deterministic covariate and $1\leq p\leq 3$;
see also \citet{CGS15} for a refinement in the case that  $p=2$ but  under the assumption of independent  and identically
distributed errors.
Furthermore, \citet[Theorem~1]{Dur02} proves that $E[\int_0^1 |\widetilde{f}_n(x)-f(x)|\, dx]=O(n^{-1/3})$.

However,  much less  is known about the asymptotic behavior of $\widetilde{f}_n$  in the case of multivariate regression models.
The only results concerning the estimator  ~$\widetilde{f}_n$
that we are aware are the following. \citet[Theorem~5]{HPW73} prove uniform consistency of~$\widetilde{f}_n$ in the
case~$d=2$ under the assumptions of deterministic regressors and a continuous target function~$f$.
Additionally, these authors  provide intuition for the convergence of
large deviation probabilities between  the estimator and  the true regression function
towards zero; see \citet[Eq.~(26)]{HPW73}.
\citet[Theorems~2.1 and~2.2]{RW75} state pointwise consistency for $\widetilde{f}_n$ in the context of a
general partial order for the regressors. Finally, \citet[Theorem~1]{CT87} prove  a consistency result
in the $d$-dimensional case. The authors  assume that the errors form a martingale difference sequence and
the covariates are  continuous and stochastic.

To the best of our knowledge, there are no any available  results concerning the integrated $L^{1}$ convergence rate of isotonic LSE in the
case of multivariate regression models.  We conjecture that a serious
obstacle for deriving such rates of convergence for  $\widetilde{f}_n$ when $d\geq 2$ is the enormous
amount of  possible lower and upper set involved in computing   (\ref{2.1}) and (\ref{2.2}); see e.g.~\citet{GW07} as well as
the discussion in Section~3 in \citet{WMO15}.

Our goal is  to fill this gap  by proposing  a suitable modification of isotonic LSE  as described in Section~\ref{S2}.
For the case of univariate regression  we let intact the isotonic LSE ~~$\widetilde{f}_n$.
However, in the multivariate case we propose a slightly simpler estimator by restricting  attention to
lower and upper sets of (hyper-)rectangular type.  As we will show, for both cases of independent regressors (see Theorem ~\ref{T2.1}) and
dependent data  (see  Theorem ~\ref{T3.1}),
such  modification avoids  the entropy problem and allows derivation of the desired convergence rate.
In sharp contrast to usual nonparametric estimators and in accordance with the classical isotonic LSE,
this estimator does not require the choice of
an appropriate bandwidth which could cause problems in our general setting with a possibly irregular
distribution of the explanatory variables and with dependent observations.
This general framework we consider  allows inclusion of
a trend component. Such a covariate  accommodates  the case of gradual changes over time in contrast to change-point models
with stationarity between these points of (abrupt) changes.

The paper is structured as follows. We introduce the proposed  estimators and present results
on their rate of convergence in Sections~\ref{S2} (independence case without trend component) and~\ref{S3}
(dependence case allowing for a deterministic trend).
Numerical examples are discussed in Section~\ref{Sec:simulations}.
All proofs as well as technical auxiliary results are deferred to Section~\ref{S4}.
\bigskip

\section{Multivariate isotonic regression under independence}
\label{S2}

Recall \eqref{eq.motivation}
where we now assume that  $f\colon \; [0,1]^d\rightarrow\R$ and
$(I_1',\varepsilon_1)',\ldots,(I_n',\varepsilon_n)'$ are independent and
identically distributed random variables on a probability space $(\Omega,{\mathcal A},P)$.
We assume that the conditional mean function~$f$ is isotonic, that is, monotonically  non-decreasing
in each argument.
Following the discussion of Section~\ref{S1}, consider the  estimators defined by
\begin{displaymath}
\underline{f}_n(x) \,\in\,
\max_{a\colon \, a\preceq x} \;\; \min_{b\colon \, b\succeq x} \mbox{Av}_Y( \bm{[} a, b \bm{]} )
\end{displaymath}
and
\begin{displaymath}
\overline{f}_n(x) \,\in\,
\min_{b\colon \, b\succeq x} \;\; \max_{a\colon \, a\preceq x} \mbox{Av}_Y( \bm{[} a, b \bm{]} ),
\end{displaymath}
where, for $a,b\in [0,1]^d$, $\bm{[} a, b \bm{]}=[a_1,b_1]\times\cdots\times [a_d,b_d]$.
It follows from the construction of both  $\underline{f}_n$ and $\overline{f}_n$ that they  are isotonic and that
$\underline{f}_n(x)\leq \overline{f}_n(x)$ holds for all~$x$.
We define the isotonic estimator ~$\widehat{f}_n$ of $f$ as any isotonic function that satisfies
\begin{equation}
\label{2.3}
\underline{f}_n(x) \,\leq\, \widehat{f}_n(x) \,\leq\, \overline{f}_n(x) \qquad \forall ~~x.
\end{equation}
For example, choose $\widehat{f}_n(x)=(\underline{f}_n(x)+\overline{f}_n(x))/2$.
In the univariate case  any choice of $\widehat{f}_n$ which is between $\underline f_n$ and $\overline f_n$
is equal to $\widetilde{f}_n$ at the observation points.
The proposed estimator  deviates from~$\widetilde{f}_n$ in the multivariate case though.
The proofs of Theorems~\ref{T2.1} and~\ref{T3.1} below show that replacing lower and upper sets by hyperrectangles  in (\ref{2.1}) and
(\ref{2.2})   simplifies the derivation of the desired rate of convergence and its computation.

Firstly, we study the case of independent and identically distributed variables $(I_t',\varepsilon_t)'$.
We impose the following condition.
\medskip
\begin{itemize}
\item[{\bf (A1)}]
\begin{itemize}
\item[(i)] The information variables $I_t$ possess a density~$p$ on~$[0,1]^d$, such that
\begin{displaymath}
C_1 \,:=\, \inf_{x\in [0,1]^d} p(x) \,\leq\, \sup_{x\in [0,1]^d} p(x) \,=:\, C_2,
\end{displaymath}
where $0<C_1\leq C_2<\infty$.
\item[(ii)] The error sequence  $(\varepsilon_t)_{t \in {\N}}$ satisfies
\begin{equation*}
E\left( \varepsilon_t \mid I_t \right)  =  0 ~~~ \mbox{a.s.}, \qquad \mbox{and} \qquad
E\left( \varepsilon_t^2 \mid I_t \right) \leq  {\overline{\sigma}_\varepsilon}^2 ~~~ \mbox{a.s.},
\end{equation*}
where ${\overline{\sigma}_\varepsilon}^2<\infty$.
\end{itemize}
\end{itemize}
\bigskip
Define
\begin{equation}
M_n \,=\, [n^{1/(d+2)}] \qquad \mbox{ and } \qquad h_n \,=\, 1/M_n.
\label{eq:bandwidth}
\end{equation}
Note that~$h_n$ corresponds to an asymptotically mean square error-optimal bandwidth of a kernel
estimator when the function to be estimated has a degree of smoothness~1. Here, and in the following, the notation  $\mathbf{1}$ denotes the $d$--dimensional vector consisting of ones.
The estimator~$\widehat{f}_n$ is based on means over hyperrectangles.
For multi-indexes $k=(k_1,\ldots,k_d)$, we define grid points by
\begin{displaymath}
x_k \,=\, (k_1 h_n,\ldots, k_d h_n)', \qquad (0\leq k_i\leq M_n \quad \forall i)
\end{displaymath}
and subsets of $[0,1]^d$ by
\begin{displaymath}
B_k \,=\, \bm{(} x_{k-\mathbf{1}}, x_k \bm{]}
\,=\, ((k_1-1)h_n, k_1h_n] \times\cdots\times ((k_d-1)h_n, k_dh_n] \qquad \forall k\in K_n,
\end{displaymath}
where $K_n=\{k\colon\; 1\leq k_1,\ldots,k_d\leq M_n\}$.
The estimator~$\widehat{f}_n$ is based on means over hyperrectangles. We expect therefore a \emph{regular} behavior
of ~$\widehat{f}_n$  if there are sufficiently many observations in each box $B_k$.
Recall that $C_1$ is the lower bound on the density of the information variables~$I_t$ which is assumed to exist by (A1)(i).
Then, regularity of ~$\widehat{f}_n$  is guaranteed to hold, provided that  the event
\begin{eqnarray}
\label{eq.an}
A_n \,=\, \left\{ \omega\colon\quad \#\{t\leq n\colon\; I_t(\omega)\in B_k\}\geq (C_1/2)\,n^{2/(d+2)} \quad
\forall k\in K_n \right\}
\end{eqnarray}
occurs with probability tending to one. This is stated in the following lemma.

{\lem
\label{L2.1}
Suppose that Assumption (A1) holds true. Then
\begin{eqnarray*}
P(A_n) \,\mathop{\longrightarrow}\limits_{n\to\infty}\, 1.
\end{eqnarray*}}
\medskip

It is well known that the traditional isotonic estimator $\widetilde{f}_n(x)$
is problematic when~$x$ is close to the boundary of the support of the distribution of the~$I_t$;
see e.g.~the discussion in \citet{SSW03}.
The same is true for~$\widehat{f}_n$ at points~$x$ near the boundary of the domain.
To fix the bias problem at extreme small and large design points,
\citet{WMO15} proposed an adequate modification by pulling up and down the isotonic LSE at these particular locations.
We do not implement a boundary correction
since this would involve some sort of tuning parameter whose appropriate choice is somewhat subjective.
In fact, we neglect the behavior of~$\widehat{f}_n$ near the boundary and focus on estimating~$f$ on the box
\begin{displaymath}
D_n \,=\, [h_n, 1-h_n]^d \,=\, \overline{ \bigcup_{k\colon\; 1<k_1,\ldots,k_d<M_n} B_k } .
\end{displaymath}
Under minimal assumptions and assuming existence of second moments for the error terms we prove
the following theorem.

{\thm
\label{T2.1}
Suppose that Assumption (A1) holds true.
Then, with $\lambda^d$ denoting the Lebesgue measure on $\R^d$,
\begin{displaymath}
E\left[ \int_{D_n} |\widehat{f}_n(x) \,-\, f(x)| \, \lambda^d(dx) \;\; \1_{A_n} \right]
\,=\, O(n^{-1/(d+2)}).
\end{displaymath}
}

First, we  notice  that  Theorem \ref{T2.1} and Lemma~\ref{L2.1} imply that
\begin{displaymath}
\int_{D_n} |\widehat{f}_n(x) \,-\, f(x)| \, \lambda^d(dx) \,=\, O_P(n^{-1/(d+2)}).
\end{displaymath}
Furthermore, consider the special case of a partially differentiable function $f\colon\, [0,1]^d\rightarrow \R$.
Then the assumption of  isotonicity implies that
$$\int_{[0,1]^d} \sum_{i=1}^d |\partial_i f(x)|\, \lambda^d(dx) \leq d(\sup_x\{f(x)\}-\inf_x\{f(x)\}).$$
Hence, the degree of smoothness, say  $\beta$, measured in the $L^1$-norm, is equal to 1.
It is well known that, under appropriate conditions, the optimal rate of convergence for the $L_1$-loss
is $n^{-\beta/(2\beta+d)}$ which reduces to $n^{-1/(d+2)}$, when $\beta=1$; see~\citet{St82}.
Hence, Theorem \ref{T2.1} indicates that $\widehat{f}_n$ achieves the
optimal rate of convergence in the class of isotonic functions.
Recall that in contrast to the classical isotonic estimator which is obtained by using
all possible lower and upper sets in (\ref{2.1}) and (\ref{2.2}) our estimator~$\widehat{f}_n$ is based
on averages over hyperrectangles only. This reduced complexity allows us to derive the desired rate of
convergence.

\bigskip

\section{Multivariate isotonic regression under dependence}
\label{S3}

Recall again \eqref{eq.motivation} where we now allow the  random variables to be dependent.
We assume the information variables to be of the form $I_{n,t}=(X_{n,t}', Z_{n,t}')'$,
where $X_{n,t}$ is a $ d_1$-dimensional vector consisting of components with values in $\N_0=\{0,1,\ldots\}$,
and $Z_{n,t}$ is a $d_2$-dimensional covariate consisting of variables with continuous marginal distribution functions and possibly a trend component~$t/n$.
Here, we allow for $d_1,\,d_2 \in \N_0$ with $d=d_1+d_2>0$.
Note that by setting $d_1=0$ or $d_2=0$, it is possible that $I_{n,t}$ is just equal to $Z_{n,t}$ or $X_{n,t}$, respectively.
More specifically, we distinguish between two cases: either the covariate vector $Z_{n,t}$ includes a trend component of the
form $t/n$, i.e. $Z_t=Z_{n,t}=(\widetilde{Z}_t',t/n)'$,
where $\widetilde{Z}_{t}$ denotes the rest of the covariates, or the covariate vector is free
of a trend.
In this section, we  consider again the isotonic estimator $\widehat{f}_{n}$ defined by \eqref{2.3}.
We show that the results of Lemma~\ref{L2.1} and  Theorem~\ref{T2.1} can be  generalized to the case of strong mixing
random variables provided that we impose some additional assumptions. We suppose that:
\begin{itemize}
\item[{\bf (A2)}]
\begin{itemize}
\item[(i)] The error sequence  $(\varepsilon_{n,t})_{t\in\N}$ satisfies
\begin{eqnarray*}
E\left( \varepsilon_{n,t} \mid I_{n,1},\ldots,I_{n,t},\varepsilon_{n,1},\ldots,\varepsilon_{n,t-1} \right)
& = & 0 \qquad \mbox{a.s.}, \\
E\left( \varepsilon_{n,t}^2 \mid I_{n,1},\ldots,I_{n,t},\varepsilon_{n,1},\ldots,\varepsilon_{n,t-1} \right)
& \leq & {\overline{\sigma}_\varepsilon}^2 \qquad \mbox{a.s.},
\end{eqnarray*}
where ${\overline{\sigma}_\varepsilon}^2<\infty$.
\item[(ii)] The process $(I_{n,t})_{t\in\N}$ is strong ($\alpha$-) mixing with corresponding mixing coefficients satisfying
\begin{displaymath}
\sum_{r=1}^\infty r^{d_2-1}\; \alpha(r) \,<\, \infty.
\end{displaymath}
\item[(iii)]
The function $f\colon D\to\R,\; D\subseteq\R^d,$ is bounded.
\end{itemize}
\end{itemize}

Having in mind that $I_{n,t}=(X_{n,t}',Z_{n,t}')'$ contains $d_1\geq 0$ components
with a discrete distribution and $d_2\geq 0$ components having either a continuous distribution or being nonrandom
such as $t/n$ we impose the following condition:\\
\begin{itemize}
\item[{\bf (A3)}]
For $t=1,\dots, n$ and $n\in \N$, the random vectors $Z_{n,t}$ consist of components with continuous marginal
distribution functions and/or a trend component $t/n$.
\begin{itemize}
\item[(i)] There exist continuous distribution functions $G_1,\ldots,G_{d_2}$ on~$\R$ and, for all $K\in\N$,
constants $0<C_1=C_1(K)\leq C_2=C_2(K)<\infty$ such that $\forall k_1,\ldots,k_{d_1}\leq K, \forall a_i\leq b_i$
\begin{eqnarray*}
\lefteqn{ C_1\, \prod_{i=1}^{d_2} \left( G_i(b_i)-G_i(a_i) \right) \quad - \quad \frac{1}{n} } \\
& \leq & \frac{1}{n} \sum_{t=1}^n P\left( X_{n,t}=(k_1, \ldots,k_{d_1})',
Z_{n,t}\in (a_1,b_1] \times \cdots \times (a_{d_2},b_{d_2}] \right) \\
& \leq & \frac{1}{n} \sum_{t=1}^n P\left( Z_{n,t}\in (a_1,b_1] \times \cdots \times (a_{d_2},b_{d_2}] \right) \\
& \leq & C_2\, \prod_{i=1}^{d_2} \left( G_i(b_i)-G_i(a_i) \right) \quad + \quad \frac{1}{n}.
\end{eqnarray*}
\item[(ii)] There exists some constant $C_3<\infty$ such that, for all $d$-dimensional
hyperrectangles~$C$,
\begin{eqnarray*}
P\Big( I_{n,t}\in C \Big| I_{n,1},\ldots,I_{n,t-d},\varepsilon_{n,1},\ldots,\varepsilon_{n,t-d} \Big)
\,\leq\, C_3 \; P\left( I_{n,t}\in C \right).
\end{eqnarray*}
\end{itemize}
\end{itemize}
\medskip

Before we proceed some comments on assumption~(A3) are in order.
Condition (A3)(i) means that the ``average distribution'' of the continuous random variables
behaves as a $d_2$-dimensional product distribution which has, after an appropriate
rescaling with $G_1^{-1},\ldots,G_{d_2}^{-1}$, a density bounded away from zero
on $[0,1]^{d_2}$. The terms $\pm 1/n$ are needed to accommodate the possible case of
a trend variable~$t/n$. Also note that assumption (A1)(i) implies the validity of assumption (A3)(ii). We impose a condition on
$P\left( I_{n,t}\in C \mid I_{n,1},\ldots,I_{n,t-d},\varepsilon_{n,1},\ldots,\varepsilon_{n,t-d} \right)$
rather than $P\left( I_{n,t}\in C \mid I_{n,1},\ldots,I_{n,t-1},\varepsilon_{n,1},\ldots,\varepsilon_{n,t-1} \right)$
in order to accommodate the case where $I_{n,t}=(Y_{n,t-1},\ldots,Y_{n,t-d})'$.

To simplify the notation, we suppress the index~$n$ in $Y_{n,t}$, $I_{n,t}$ and $\varepsilon_{n,t}$ from here on,
just keeping in mind that also a triangular scheme is allowed, e.g., when a trend variable~$t/n$ is included.
We define
\begin{displaymath}
\widetilde{M}_n \,=\, [n^{1/(d_2+2)}], \qquad \widetilde{h}_n \,=\, 1/\widetilde{M}_n,
\end{displaymath}
$d=d_1+d_2$ and, for multi-indexes $k=(k_1,\ldots,k_d)$
($0\leq k_j\leq K$ $\forall j=1,\ldots,d_1$, $1\leq k_j\leq \widetilde{M}_n$ $\forall j=d_1+1,\ldots,d$),
subsets of the domain of~$f$ as
\begin{displaymath}
\widetilde{B}_k \,=\, \{(k_1,\ldots,k_{d_1})'\}\times (G_1^{-1}((k_{d_1+1}-1)\widetilde{h}_n),G_1^{-1}(k_{d_1+1}\widetilde{h}_n)]
\times\cdots\times (G_{d_2}^{-1}((k_d-1)\widetilde{h}_n),G_{d_2}^{-1}(k_d\widetilde{h}_n)].
\end{displaymath}

Since the estimator~$\widehat{f}_n$ is based on means over hyperrectangles, a ``regular'' behavior
of it can be expected if there are sufficiently many observations in each box $B_k$.
It turns out that regularity can be assured if the following event occurs:
\begin{eqnarray}
\label{eq.an-tilde}
\widetilde{A}_n \,=\, \left\{ \omega\colon\; \#\{t\leq n\colon\; I_t(\omega)\in \widetilde{B}_k\}\geq C_4\,n^{2/(d_2+2)} \quad
\forall k\in \widetilde{K}_n \right\},
\end{eqnarray}
where $C_4$ is some positive constant
and $\widetilde{K}_n=\{0,\ldots,K\}^{d_1}\times\{1,\ldots,\widetilde{M}_n\}^{d_2}$.

Using a Fuk-Nagaev-type inequality for dependent random variables we can prove the
following analogous result to Lemma~\ref{L2.1}.

{\lem
\label{L3.1}
Suppose that Assumptions (A2) and (A3) hold true. Then, for sufficiently small $C_4>0$ in (\ref{eq.an-tilde}),
\begin{eqnarray*}
P(\widetilde{A}_n) \,\mathop{\longrightarrow}\limits_{n\to\infty}\, 1.
\end{eqnarray*}
}
\medskip

Recall again that the traditional isotonic estimator $\widetilde{f}_n(x)$
is problematic when~$x$ is close to the boundary of the support of the distribution of the~$I_t$.
We neglect the behavior of~$\widehat{f}_n$ near the boundary and focus on estimating~$f$ on
\begin{displaymath}
\widetilde{D}_n
\,=\, \{0,\ldots,K\}^{d_1}\times (G_1^{-1}(\widetilde{h}_n),G_1^{-1}(1-\widetilde{h}_n)]
\times\cdots\times(G_{d_2}^{-1}(\widetilde{h}_n),G_{d_2}^{-1}(1-\widetilde{h}_n)].
\end{displaymath}
Denote by $Q_1,\ldots,Q_{d_2}$ the probability measures corresponding to the
distribution functions $G_1,\ldots,G_{d_2}$, respectively. With $\mu^{d_1}$ being the
counting measure on~$\N_0^{d_1}$, define $\nu=\mu^{d_1}\otimes Q_1\otimes\cdots\otimes Q_{d_2}$.

{\thm
\label{T3.1}
Suppose that Assumptions (A2) and (A3) hold true.
Then,
\begin{displaymath}
E\left[ \int_{\widetilde{D}_n} |\widehat{f}_n(x) \,-\, f(x)| \, \nu(dx) \;\; \1_{\widetilde{A}_n} \right]
\,=\, O(n^{-1/(d_2+2)}).
\end{displaymath}
Here, in the definition of the event $\widetilde{A}_n$, the constant $C_4$ is chosen such that the assertion of Lemma~\ref{L3.1} holds true.
}
\medskip

\noindent
Some remarks are in order. First, it follows again  from this theorem and Lemma~\ref{L3.1} that
\begin{displaymath}
\int_{\widetilde{D}_n} |\widehat{f}_n(x) \,-\, f(x)| \, \nu(dx) \,=\, O_P(n^{-1/(d_2+2)}).
\end{displaymath}
Furthermore, we point out that the obtained rate of convergence does not depend on the number of discrete
explanatory random variables. This is explained by the fact that, for any $k_1,\ldots,k_{d_1}\in \{0,\ldots,K\}$,
the cardinality of the set $\{t\leq n\colon\; X_{n,t}=(k_1,\ldots,k_{d_1})'\}$ is proportional to the
sample size~$n$. Therefore, there is no need to smooth over the first $d_1$ directions and there is
no loss due to a trade-off between bias and variance that would appear with nonparametric smoothing
techniques.
\medskip

Properties of the noise process can be  taken into account, provided that we  have some prior knowledge.
Indeed, if we knew  the conditional variance $E(\epsilon_t^2\mid I_t)$,
e.g.~in the case of a known distributional family for the errors, then we could replace
the  means $\mbox{Av}_Y(B) = \sum_{t=1}^n Y_t \1( I_t \in B ) / \#\{t\leq n\colon \; I_t\in B\}$
by the weighted means $\sum_{t=1}^n w(I_t) Y_t \1( I_t \in B ) / \sum_{t=1}^n w(I_t) \1( I_t \in B )$,
where the weights $w(I_t)$ are proportional to $1/E(\epsilon_t^2\mid I_t)$.
This corresponds to a weighted least squares estimator in linear regression.
However, our main intention was to produce a general, fully nonparametric method. Since prior
knowledge of $E(\epsilon_t^2\mid I_t)$ is rarely available, we pursue  the  approach
based on unweighted means.
\medskip

{\exmp
Suppose that integer-valued random variables $Y_1,\ldots,Y_n$ are observed, where
\begin{displaymath}
Y_t \mid {\mathcal F}_{t-1} \,\sim\, \mbox{Poisson}(\lambda_t),
\end{displaymath}
\begin{displaymath}
\lambda_t \,=\, f( Y_{t-1}, Z_{t-1} ),
\end{displaymath}
and $Z_t$ is a covariate with values in $[0,1]^{d_2}$ which is independent of $Y_t,\ldots,Y_0,Z_{t-1},\ldots,Z_0$,
${\mathcal F}_s=\sigma(Y_s,Z_s,\ldots,Y_0,Z_0)$.
Assume that the function $f\colon \; \N_0\times [0,1]^{d_2}\rightarrow [0,M]$ is isotonic and bounded by $M<\infty$.
The information variable at time~$t$ is $I_t=(Y_{t-1},Z_{t-1}')'$.
We have that
\begin{displaymath}
Y_t \,=\, f(I_t) \,+\, \varepsilon_t,
\end{displaymath}
where
\begin{eqnarray*}
E\left( \varepsilon_t \mid I_t,\ldots,I_0,\varepsilon_{t-1},\ldots,\varepsilon_{0} \right) & = & 0, \\
E\left( \varepsilon_t^2 \mid I_t,\ldots,I_0,\varepsilon_{t-1},\ldots,\varepsilon_{0} \right)
& = & f(I_t) \,\leq\, M.
\end{eqnarray*}
It can be shown that Assumption (A2)(i) is also fulfilled.
Indeed, let $Q_t^k := P^{Y_t\mid Y_{t-1}=k} = \int \mbox{Poisson}(f(k,z))\, P^{Z_{t-1}}(dz)$.
Since $f(k,z)\in[0,M]$ for all values of $k$ and $z$,
\begin{displaymath}
\inf_t \inf_{k\in \N_0} Q_t^k(\{0\}) \,>\, 0,
\end{displaymath}
that is, Doeblin's condition is fulfilled.
It follows from Theorem~2.4.1 on page~88 in \citet{Dou94} that the Markov chain $(Y_t)_t$ is uniformly ($\phi$-) mixing
and, therefore, absolutely regular with coefficients satisfying
\begin{displaymath}
\beta^Y(k) \,\leq\, C \rho^k \quad \forall k\in\N_0,
\end{displaymath}
for some $C<\infty$ and $\rho\in [0,1)$.
Since the process $(I_t)_t$ is also a Markov chain, we obtain that
\begin{eqnarray*}
\lefteqn{ \beta\left( \sigma(I_0,I_1,\ldots,I_t), \sigma(I_{t+k},I_{t+k+1},\ldots) \right) } \\
& = & \beta\left( \sigma(I_t), \sigma(I_{t+k}) \right) \\
& \leq & \beta\left( \sigma(Y_t,I_t), \sigma(I_{t+k}) \right) \\
& = & \beta\left( \sigma(Y_t), \sigma(I_{t+k}) \right) \\
& = & \beta\left( \sigma(Y_t), \sigma(Y_{t+k-1}) \right).
\end{eqnarray*}
(The first and the last but one equalities follows from the Markovian structure;
see also the note after Theorem~7.3 in \citet{Bra07}. The last one follows from
independence of $Z_{t+k-1}$ and $(Y_t,Y_{t+k-1})$; see also Theorem~6.2 in \citet{Bra07}.)
Hence, the coefficients of absolute regularity of the process $(I_t)_t$ satisfy
\begin{displaymath}
\beta^{I}(k) \,\leq\, \beta^Y(k-1) \,\leq\, C\rho^{k-1} \quad \forall k\in\N.
\end{displaymath}
}
\bigskip

\section{Applications}
\label{Sec:simulations}

\subsection{Simulations}
\label{Sub-Sec:sims}

\begin{figure}
    \centering
    \includegraphics[width=0.7\textwidth]{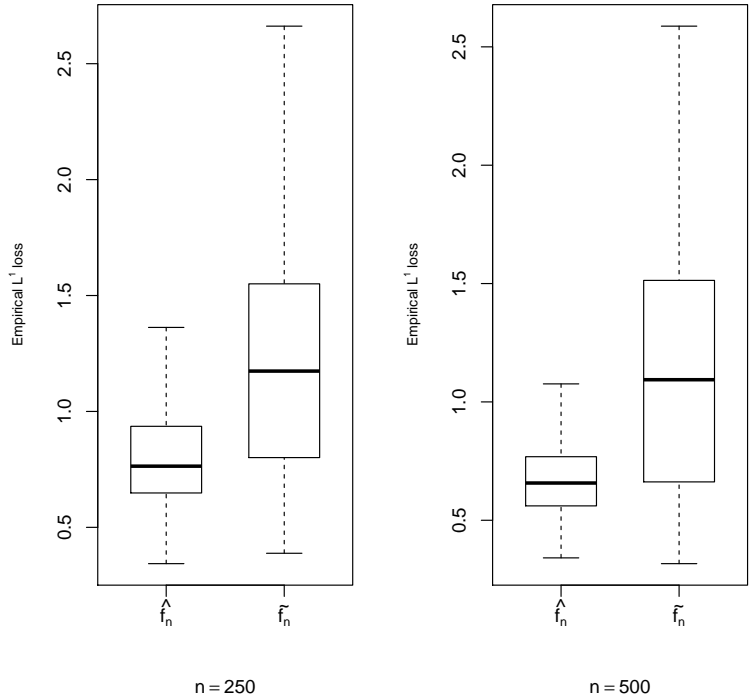}
\caption{Box plots of empirical $L^{1}$ loss values for $\widehat{f}_{n}$ and the isotonic LSE $\widetilde{f}_{n}$.}
\label{fig:Simulation}
\end{figure}

We illustrate the theoretical results by a limited simulation study comparing the performance of $\widehat{f}_{n}$ and the isotonic LSE $\widetilde{f}_{n}$ in terms of their $L^{1}$ error.
More specifically, consider the following  Poisson count time series model, as described in Example~3.1, where we have assumed that $Z_{t}$
is a deterministic trend, i.e. for $t=1,2,\ldots, n$
\begin{displaymath}
Y_t \mid {\mathcal F}_{t-1} \,\sim\, \mbox{Poisson}(\lambda_t),~~
\lambda_t \,=\, f( Y_{t-1}, Z_{t-1} ),
\end{displaymath}
where
$$
f(y,z) = -5 + \frac{20}{1+\exp(-0.3y)} + 4z\qquad\text{ and }\qquad Z_t=\frac{t}{n}.
$$
We compute the isotonic least squares estimator by using the R package \texttt{isotonic.pen}
which returns  the values of the estimated function on an equidistant $21\,\times\, 21$ grid; see \citet{WMO15} for details.
Note that $\widehat f_n=(\bar f_n+\underline f_n)/2$  can be computed (exactly) on any grid. We use the grid employed  by  \texttt{isotonic.pen}
discarding some  points so that we can  avoid any  boundary issues. To this end, we choose the lower left\;/\;upper  right corner of the grid such that, on the one hand,  the number of observed
information variables within the corresponding rectangle  is maximized  and on the other hand, for every grid point all  upper and lower rectangles contain at least one data point.
To compare the empirical  performance of the estimators, we compute the integrated  $L^{1}$ error over the grid values. This process is repeated 500 times.
Figure~\ref{fig:Simulation} shows box plots of the values of integrated $L^{1}$ error for two different sample sizes and illustrates 
that $\widehat{f}_{n}$ achieves smaller error than the isotonic LSE.

\subsection{Data Examples}
\label{Sub-Sec:data}

\begin{figure}
    \centering
    \includegraphics[width=\textwidth]{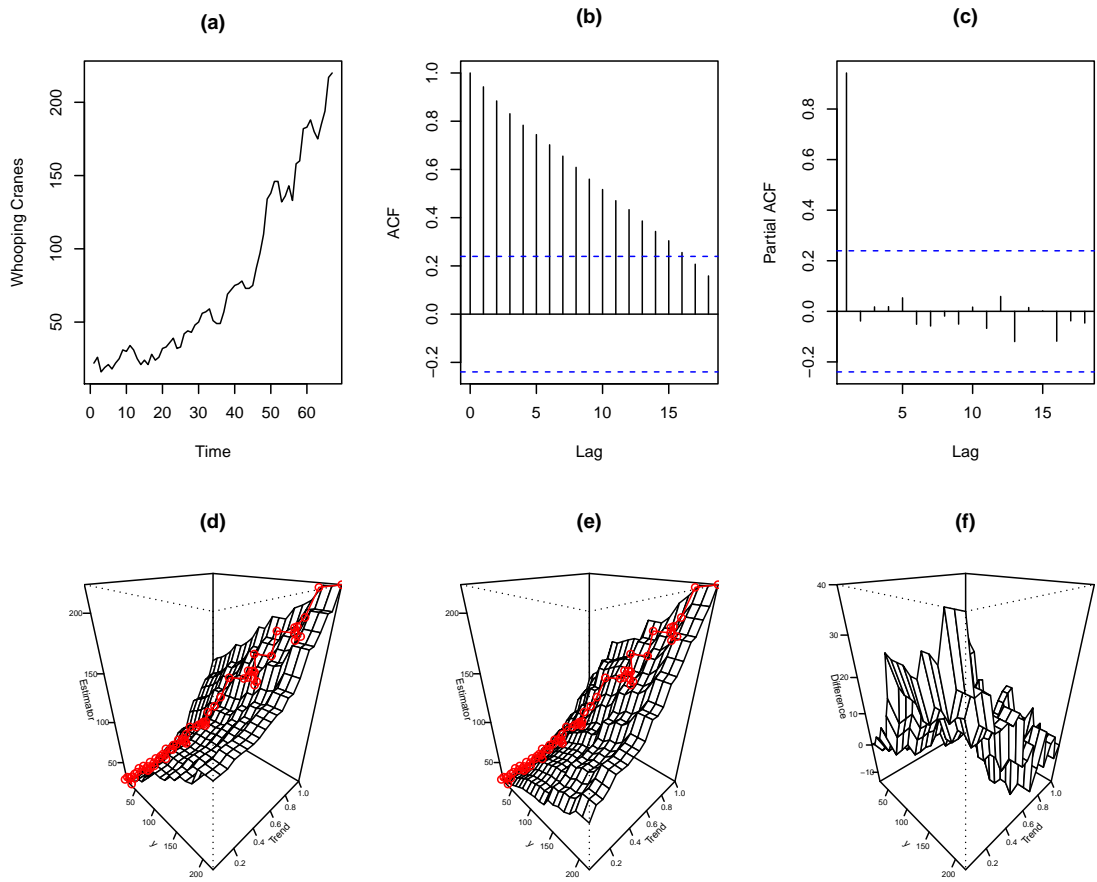}
\caption{(a) Time series plot of the yearly number of whooping cranes between 1938 to 2005. (b) Autocorrelation function. (c) Partial autocorrelation
function. (d) Plot of isotonic LSE $\widetilde{f}_{n}$ and the data (red points). (e) Plot of the estimator $\widehat{f}_{n}$ and the data (red points). (f) Plot of the difference $\widetilde{f}_{n}- \widehat{f}_{n}$. }
\label{fig:Data1}
\end{figure}

We apply the methodology  to biological and financial time series which exhibit some form of non-stationarity.
First, we investigate the population growth of whooping cranes that became nearly extinct during the period 1938-1955.
Whooping cranes are one of the largest birds in North America but also one of the rarest that can be found in the continent.
For some time their population has been constantly decreasing and reached to about 20 individual birds in the world. With
the employment of various conservation measures the population grew over the last years. The data we have are depicted in
Figure~\ref{fig:Data1} which shows the growth of population of whooping cranes between 1938 to 2005; see \citet{data}. Note that this is a case of an integer valued time series.
The second example  refers to daily net asset value (NAV) of the  BlackRock  Global Allocation Fund  during the period  1/4/2016 to 30/1/2018. Here we note that series takes
values on real numbers. 

\begin{figure}
    \centering
    \includegraphics[width=\textwidth]{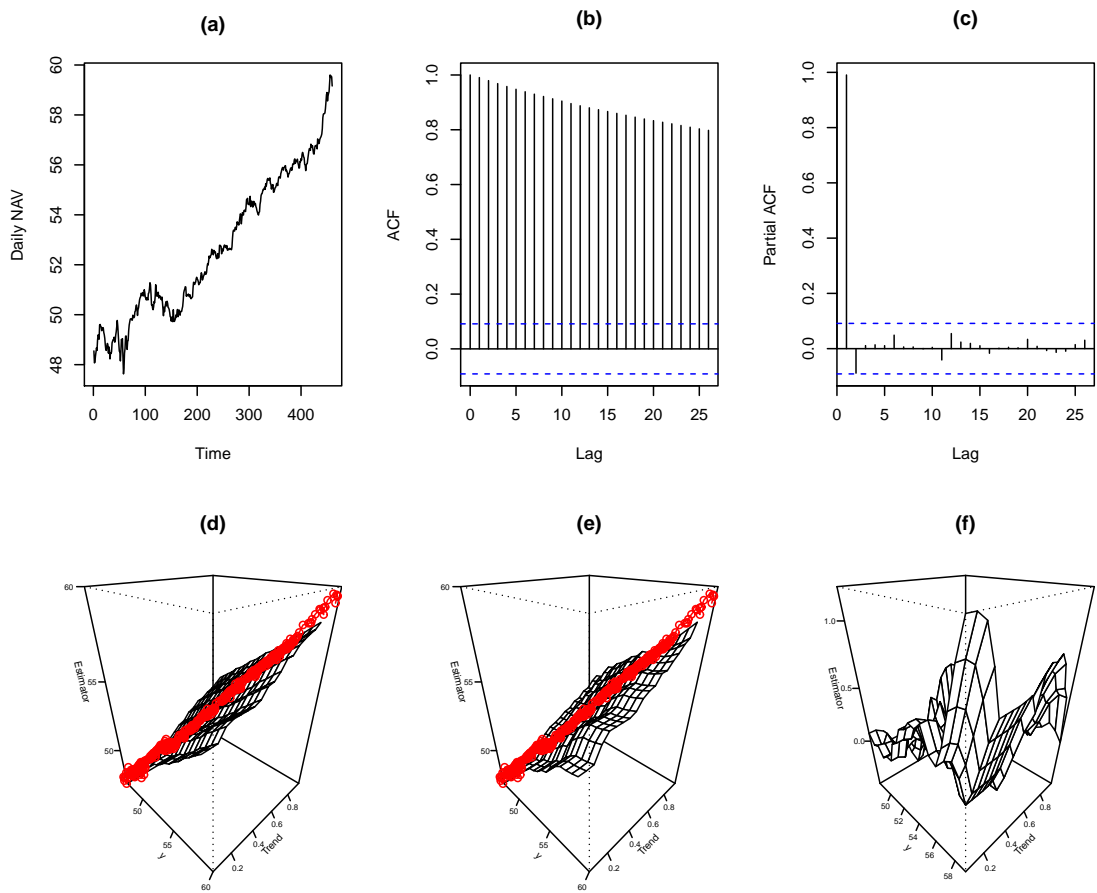}
\caption{(a) Time series plot of daily NAV prices of BlackRock  Global Allocation Fund  during the period  1/4/2016 to 30/1/2018.
(b) Autocorrelation function. (c) Partial autocorrelation
function. (d) Plot of isotonic LSE $\widetilde{f}_{n}$ and the data (red points). (e) Plot of the estimator $\widehat{f}_{n}$ and the data (red points). (f) Plot of the difference $\widetilde{f}_{n}- \widehat{f}_{n}$. }
\label{fig:Data2}
\end{figure}

For both of these data examples, a simple time series plot   reveals  increasing trend and  strong autocorrelation which decays slowly. 
The partial autocorrelation functions shows a strong autocorrelation at lag~1; see the upper panel of Figures~\ref{fig:Data1}  and~\ref{fig:Data2}.
We fit a non-parametric time series model  to these data by using isotonic estimation methods. We include the covariate vector
$I_t=(Y_{t-1}, t/n)'$, where $n$ is number of effective observations (e.g.~for the population growth of whooping cranes  the number of observation  is equal to  68 but $n=67$ because of the inclusion of $Y_{t-1}$).
We consider again the estimator  $\widehat{f}_{n}$  and the isotonic LSE $\widetilde f_n$ and work the same way as it was explained  in Subsection \ref{Sub-Sec:sims} .
The lower panels of Figures~\ref{fig:Data1} and \ref{fig:Data2} show that both estimators are quite close  near the observation points, 
but differ significantly at some grid points that are located far from the bulk of data.
We examine the performance of both methods  for estimating the two models.
This task is accomplished by studying  the in sample predictive
power using the mean absolute prediction error (MAPE), that is $\sum_{t=1}^{n} |\widehat{Y}_{t}- Y_{t}|/n$. Here, $\widehat Y_t$ is obtained by evaluating $\widehat f_n$ and $\widetilde f_n$, respectively, on a grid point close to $(Y_{t-1}, t/n)'$. The results are shown in Table \ref{Table:data}. Clearly, the new estimator $\widehat{f}_{n}$ outperforms the isotonic LSE in terms of MAPE for both data examples.

%
\begin{table}
\begin{tabular}{|l||c|c|}
  \hline
	                 & Example 1 (Whooping cranes)	 & Example 2 (BlackRock
Global Allocation Fund)\\\hline
  $\widehat f_n$   &  5.916453 &  0.3108097   \\ 
  $\widetilde f_n$ &  6.002488     &  0.3118416       \\
  \hline
\end{tabular}
\caption{MAPE obtained from our estimator $\widehat f_n$ and the isotonic LSE  $\widetilde{f}_{n}$.}
\label{Table:data}
\end{table}

\clearpage

\section{Proofs and Auxiliary Results}
\label{S4}

\noindent
We prove our main results in Section~\ref{SS4.1}.
Some auxiliary lemmas are stated and proved in Section~\ref{SS4.2}.
\medskip

\subsection{Proofs of the main results}
\label{SS4.1}

\begin{proof}[Proof of Lemma~\ref{L2.1}]
Since $\sum_{t=1}^n P(I_t \in B_k)\geq nC_1 h_n^d\geq C_1 n^{2/(d+2)}$ it suffices to show that
\begin{equation}
\label{pl21.1}
\sum_{k\in K_n} P\left( \sum_{t=1}^n P(I_t\in B_k) \,-\, \1(I_t\in B_k)
\,\geq\, \frac{C_1\; n^{2/(d+2)}}{2} \right)
\,\mathop{\longrightarrow}\limits_{n\to\infty}\, 0,
\end{equation}
by recalling that  $\1(.)$ denotes the indicator function.
We obtain from Bernstein's inequality, for all $k\in K_n$
and $\mu_n=C_1 n^{2/(d+2)}/2$, that
\begin{eqnarray*}
P\left( \sum_{t=1}^n P(I_t\in B_k) \,-\, \1(I_t\in B_k)
\,\geq\, \mu_n \right)
& \leq & \exp\left\{ - \; \frac{ \mu_n^2/2 }
{ n\; P(I_1\in B_k) \; (1-P(I_1\in B_k)) \,+\, \mu_n/3 } \right\} \\
& \leq & \exp\left\{ - \; C \; n^{2/(d+2)} \right\},
\end{eqnarray*}
for some $C>0$, which proves (\ref{pl21.1}).
\end{proof}

\begin{proof}[Proof of Theorem~\ref{T2.1}]
We analyze the contribution of the stochastic part and the bias of the estimator
separately. For the latter, we exploit the assumed isotonicity in conjunction
with boundedness of~$f$ in order to construct an  estimate of the integrated bias from above and below.
To this end, denote for an arbitrary function~$g$ its positive (respectively negative) part by $g_{+}$ (respectively $g_{-}$).
Then, it suffices to show that
\begin{subequations}
\begin{align}
\label{pt21.1}
E\left[ \int_{D_n} \left( \widehat{f}_n(x) \,-\, f(x) \right)_+ \, \lambda^d(dx) \;\; \1_{A_n} \right]
& =  O\left( n^{-1/(d+2)} \right),
\end{align}
and
\begin{align}
\label{pt21.2}
E\left[ \int_{D_n} \left( \widehat{f}_n(x) \,-\, f(x) \right)_- \, \lambda^d(dx) \;\; \1_{A_n} \right]
& =  O\left( n^{-1/(d+2)} \right).
\end{align}
\end{subequations}
We have, for all $x\in B_k=\bm{(} x_{k-\mathbf{1}},x_k\bm{]}$,
\begin{eqnarray}
\label{pt21.3}
\left( \widehat{f}_n(x) \,-\, f(x) \right)_+
& \leq & \left( \sup_{y\preceq x_k} \mbox{Av}_Y( \bm{(} y, x_{k+\mathbf{1}} \bm{]} ) \,-\, f(x) \right)_+
\nonumber \\
& \leq & \sup_{y\preceq x_k} \left| \mbox{Av}_\varepsilon( \bm{(} y, x_{k+\mathbf{1}} \bm{]} ) \right|
\,+\, \left( f(x_{k+\mathbf{1}}) \,-\, f(x) \right).
\end{eqnarray}
By Lemma~\ref{LA.1} and since $\lambda^d(B_k)=h_n^d=O(n^{-d/(d+2)})$ we obtain for the bias that
\begin{eqnarray}
\label{pt21.4}
\lefteqn{ \sum_{k\colon \, 1<k_1,\ldots,k_d<M_n}
\int_{B_k} \left( f(x_{k+\mathbf{1}}) \,-\, f(x) \right) \, \lambda^d(dx) } \nonumber \\
& \leq & \sum_{k\colon \, 1<k_1,\ldots,k_d<M_n} \left( f(x_{k+\mathbf{1}}) \,-\, f(x_{k-\mathbf{1}}) \right) \; \lambda^d(B_k)
\,=\, O\left( n^{-1/(d+2)} \right).
\end{eqnarray}

\noindent
For the stochastic part, we estimate $\; E[\sup_{y\preceq x_k} |\mbox{Av}_\varepsilon( \bm{(} y, x_{k+\mathbf{1}} \bm{]} )|]$.
For this purpose, we define a dyadic scheme of nested hyperrectangles: For $j_1,\ldots,j_d\geq 0$,
\begin{displaymath}
B_k^{(j_1,\ldots,j_d)}
\,=\, ((k_1+1-2^{j_1})h_n, (k_1+1)h_n] \times \cdots \times ((k_d+1-2^{j_d})h_n, (k_d+1)h_n].
\end{displaymath}
(We have in particular $B_k^{(0,\ldots,0)}=B_{k+\mathbf{1}}$ and
$B_k^{(j_1,\ldots,j_d)}= \bigcup_{m_1,\ldots,m_d\colon\, 0\leq m_i\leq 2^{j_i}-1}
B_{(k_1+1-m_1,\ldots,k_d+1-m_d)}$.)
Since $\{y\preceq x_k\}\subseteq \bigcup_{j_1,\ldots,j_d\geq 0} B_k^{(j_1+1,\ldots,j_d+1)}\setminus B_k^{(j_1,\ldots,j_d)}$
and since $y\in B_k^{(j_1+1,\ldots,j_d+1)}\setminus B_k^{(j_1,\ldots,j_d)}$ implies that
$\bm{(} y, x_{k+\mathbf{1}} \bm{]}\supseteq B_k^{(j_1,\ldots,j_d)}$
we obtain, for all $x\in B_k$,
\begin{displaymath}
\sup_{y\preceq x_k} \left| \mbox{Av}_\varepsilon( \bm{(} y, x_{k+\mathbf{1}} \bm{]} ) \right|
\,\leq\, \sum_{j_1,\ldots,j_d\geq 0}
\frac{ \sup\left\{ \left| \sum_{t=1}^n \varepsilon_t \1( I_t\in \bm{(} y,x_{k+\mathbf{1}} \bm{]} )
\right|\colon \;\; y\in B_k^{(j_1+1,\ldots,j_d+1)}\setminus B_k^{(j_1,\ldots,j_d)} \right\} }
{ \# \left\{ t\leq n\colon \;\; I_t\in B_k^{(j_1,\ldots,j_d)} \right\} }.
\end{displaymath}
Recall that if the event~$A_n$ occurs, then $\#\{t\leq n\colon \; I_t\in B_k\}\geq (C_1/2)\, n^{2/(d+2)}$
for all~$k\in K_n$, which implies that
\begin{displaymath}
\# \left\{t\leq n\colon \;\; I_t\in B_k^{(j_1,\ldots,j_d)} \right\} \,\geq\, (C_1/2) \; 2^{j_1+\cdots +j_d}
\; n^{2/(d+2)}.
\end{displaymath}
Furthermore, it follows from Lemma~\ref{LA.2} that for some $C<\infty$
\begin{equation}
\label{pt21.5}
E\left[ \sup\left\{ \left| \sum_{t=1}^n \varepsilon_t \1( I_t\in \bm{(} y,x_{k+\mathbf{1}} \bm{]} )
\right|\colon \;\; y\in B_k^{(j_1+1,\ldots,j_d+1)} \right\} \right]
\,\leq\, C \; 2^{(j_1+\cdots +j_d)/2} \; n^{1/(d+2)}.
\end{equation}
Therefore, we obtain that
\begin{eqnarray*}
\lefteqn{ \sup_{x\in B_k} E\left[ \sup_{y\preceq x_k}
\left| \mbox{Av}_\varepsilon( \bm{(} y, x_{k+\mathbf{1}} \bm{]} ) \right| \;\; \1_{A_n} \right] } \nonumber \\
& \leq & \sum_{j_1,\ldots,j_d\geq 0} E\left[
\frac{ \sup \left\{ \left| \sum_{t=1}^n \varepsilon_t \; \1(I_t\in \bm{(} y, x_{k+\mathbf{1}} \bm{]} ) \right|
\colon \quad y\in B_k^{(j_1+1,\ldots,j_d+1)} \right\} }
{ \#\{t\leq n\colon\; I_t\in B_k^{(j_1,\ldots,j_d)}\} }
\;\; \1_{A_n} \right] \nonumber \\
& = & O\left( n^{-1/(d+2)} \; \sum_{j_1,\ldots,j_d\geq 0} 2^{-(j_1+\cdots +j_d)/2} \right)
\,=\, O\left( n^{-1/(d+2)} \right).
\end{eqnarray*}
This yields, in conjunction with (\ref{pt21.3}) and (\ref{pt21.4}), that (\ref{pt21.1}) holds.
The proof of (\ref{pt21.2}) is completely analogous and therefore it is  omitted.
\end{proof}
\medskip

\begin{proof}[Proof of Lemma~\ref{L3.1}]
Analogously to the proof of Lemma~\ref{L2.1}, we will show that
\begin{equation}
\label{pl31.1}
\max_{k\in \widetilde{K}_n} \left\{ P\left( \sum_{t=1}^n P(I_t\in \widetilde{B}_k) \,-\, \1(I_t\in \widetilde{B}_k)
\,\geq\, \frac{C_1\; n^{2/(d_2+2)}}{2} \right) \right\}
\,=\, o\left( n^{-d_2/(d_2+2)} \right).
\end{equation}
Let $\mu_n= {C_1 n^{2/(d_2+2)}}/{8}$ and, for arbitrary $k\in \widetilde{K}_n$,
$\eta_t=\1(I_t\in \widetilde{B}_k)-P(I_t\in \widetilde{B}_k)$.
It follows from the Fuk-Nagaev-type inequality (I.6) of \citet[page~4]{Rio00} that,
for all $\kappa\geq 1$,
\begin{equation}
\label{pl31.2}
P\left( \left| \sum_{t=1}^n \eta_t \right| \,\geq\, 4\, \mu_n \right)
\,\leq\, \; \left( \left( 1 \,+\, \frac{\mu_n^2}{\kappa\; s_n^2} \right)^{-\kappa/2}
\,+\, \frac{ n\; \alpha([\mu_n/\kappa]) }{ \mu_n } \right),
\end{equation}
where
\begin{displaymath}
s_n^2 \,=\, \sum_{s,t=1}^n \left| \cov( \eta_s, \eta_t ) \right|.
\end{displaymath}

Note that (A2)(ii) implies that
\begin{equation}\label{eq.alpha}
\alpha(N)=o(N^{-d_2})\qquad\text{and}\qquad\sum_{r=N}^\infty \alpha(r)=o(N^{1-d_2}).
\end{equation}
To see the first relationship of \eqref{eq.alpha}, note that, for $j\in\N$,
\begin{eqnarray*}
\alpha( 2^j )
& \leq & \frac{1}{2^{j-1}} \sum_{r=2^{j-1}+1}^{2^j} \alpha(r)
\leq  \frac{1}{2^{j-1}} \sum_{r=2^{j-1}+1}^{2^j} \left( \frac{r}{2^{j-1}} \right)^{d_2-1} \alpha(r) \\
& = & 2^{(1-j)d_2} \sum_{r=2^{j-1}+1}^{2^j} r^{d_2-1} \alpha(r),
\end{eqnarray*}
which implies that $2^{d_2j}\alpha(2^j)\mathop{\longrightarrow}_{j\to\infty} 0$.
Since the sequence $(\alpha(r))_{r\in\N}$ is monotonically non-increasing, we obtain that $\alpha(N)=o(N^{-d_2})$.
Note that the second relationship of \eqref{eq.alpha} follows from
$N^{d_2-1}\sum_{r=N}^\infty \alpha(r) \,\leq\, \sum_{r=N}^\infty r^{d_2-1}\alpha(r) \mathop{\longrightarrow}_{j\to\infty} 0$.

Since $\alpha([\mu_n/\kappa])=o(n^{-2d_2/(d_2+2)})$ we obtain that
\begin{displaymath}
\frac{ n\; \alpha([\mu_n/\kappa]) }{ \mu_n }
\,=\, o\left( n\; n^{-2d_2/(d_2+2)}\; n^{-2/(d_2+2)} \right) \,=\, o\left( n^{-d_2/(d_2+2)} \right),
\end{displaymath}
that is, the second term on the right-hand side of (\ref{pl31.2}) is of the required order.
It remains to estimate $s_n^2$. To this end, we distinguish between the two cases of covariates  without
and with a trend component.
In the first case, we obtain from the upper bound in (A3)(i) that, for all $t,r$ with $1\leq t\leq t+r\leq n$,
\begin{eqnarray*}
\left| \cov( \eta_t, \eta_{t+r} ) \right| \,=\, \left\{ \begin{array}{ll}
O\left( n^{-d_2/(d_2+2)} \right) & \quad \mbox{ if } 0\leq r<d, \\
O\left( n^{-2d_2/(d_2+2)} \right) & \quad \mbox{ if } r\geq d
\end{array} \right. .
\end{eqnarray*}
On the other hand, we obtain from a covariance inequality for strong mixing processes
(see e.g.~\citet[Corollary~10.16]{Bra07}) that
\begin{displaymath}
\left| \cov( \eta_t, \eta_{t+r} ) \right| \,\leq\, 4\; \alpha(r)\; \|\eta_t\|_\infty
\; \|\eta_{t+r}\|_\infty \,\leq\, 4\; \alpha(r).
\end{displaymath}
Therefore, with $N_n=[n^{d_2/(d_2+2)}]$,
\begin{eqnarray}
\label{pl31.5}
s_n^2
& \leq &
\sum_{s,t\colon\; |s-t|<d} \left| \cov( \eta_s, \eta_t ) \right|
\,+\, \sum_{s,t\colon\; d\leq |s-t|\leq N_n} \left| \cov( \eta_s, \eta_t ) \right|
\,+\, 2n\; \sum_{r=N_n+1}^{n-1} 4 \; \alpha(r)  \nonumber \\
& = & O\left( n^{2/(d_2+2)} \right) \,+\, O\left( n^{2/(d_2+2)} \right) \,+\, o\left( n\; N_n^{1-d_2} \right) \nonumber \\
& = & O\left( n^{2/(d_2+2)} \right) \,+\, o\left( n^{2/(d_2+2)} \; n^{(2-d_2)d_2/(d_2+2)} \right).
\end{eqnarray}
In the case with trend, we get from (A3)(i) that
\begin{eqnarray*}
\left| \cov( \eta_t, \eta_{t+r} ) \right| \,=\, \left\{ \begin{array}{ll}
O\left( n^{-(d_2-1)/(d_2+2)} \right) & \quad \mbox{ if } 0\leq r<d-1, \\
O\left( n^{-2(d_2-1)/(d_2+2)} \right) & \quad \mbox{ if } r\geq d-1
\end{array} \right. .
\end{eqnarray*}
On the other hand, we see that $I_t=(\widetilde{I}_t',t/n)'\not\in \widetilde{B}_k$,
and therefore $\eta_t=0$ if $t\not\in I_{n,k}:= ((k_d-1)nh_n,k_dnh_n]$.
Hence, here with $N_n=[n^{(d_2-1)/(d_2+2)}]$,
\begin{eqnarray}
\label{pl31.6}
s_n^2
& \leq & \sum_{s,t\in I_{n,k}} \left| \cov( \eta_s, \eta_t ) \right| \nonumber \\
& = & O\left( \sum_{(s,t)\in I_{n,k}\colon\; |s-t|<d-1} \left| \cov( \eta_s, \eta_t ) \right|
\,+\, \sum_{(s,t)\in I_{n,k}\colon\; d-1\leq |s-t|\leq N_n} \left| \cov( \eta_s, \eta_t ) \right|
\,+\, nh_n\; \sum_{r=N_n+1}^{n-1} 4 \; \alpha(r) \right) \nonumber \\
& = & O\left( n^{2/(d_2+2)} \right) \,+\, o\left( n^{2/(d_2+2)} \; n^{(2-d_2)(d_2-1)/(d_2+2)} \right).
\end{eqnarray}
We see from (\ref{pl31.5}) and (\ref{pl31.6}) that in the two cases without and with trend the term
$(1 + \mu_n^2/(\kappa s_n^2))^{-1}$ is of order $O(n^{-\gamma})$, for some $\gamma>0$.
Choosing $\kappa>2d_2/\gamma$ we see that (\ref{pl31.1}) follows from (\ref{pl31.2}), which completes the proof.
\end{proof}
\medskip

\begin{proof}[Proof of Theorem~\ref{T3.1}]
The proof of this theorem is largely the same as that of Theorem~\ref{T2.1}.
We show that
\begin{subequations}
\begin{equation}
\label{pt31.1}
E\left[ \int_{\widetilde{D}_n} \left( \widehat{f}_n(x) \,-\, f(x) \right)_+ \, \nu(dx) \;\; \1_{\widetilde{A}_n} \right]
\,=\,  O\left( n^{-1/(d_2+2)} \right),
\end{equation}
and
\begin{equation}
\label{pt31.2}
E\left[ \int_{\widetilde{D}_n} \left( \widehat{f}_n(x) \,-\, f(x) \right)_- \, \nu(dx) \;\; \1_{\widetilde{A}_n} \right]
\,=\, O\left( n^{-1/(d_2+2)} \right).
\end{equation}
\end{subequations}
We define grid points
\begin{eqnarray*}
x_k & = & \left( k_1,\ldots,k_{d_1},G_1^{-1}(k_{d_1+1} h_n),\ldots,G_{d_2}^{-1}(k_d h_n) \right)', \\
\overline{x}_k
& = & \left( k_1,\ldots,k_{d_1},G_1^{-1}((k_{d_1+1}+1) h_n),\ldots,G_{d_2}^{-1}((k_d+1) h_n) \right)', \\
\underline{x}_k
& = & \left( k_1,\ldots,k_{d_1},G_1^{-1}((k_{d_1+1}-1) h_n),\ldots,G_{d_2}^{-1}((k_d-1) h_n) \right)'. \\
\end{eqnarray*}
We have, for all $x\in \widetilde{B}_k=\bm{(} \underline{x}_k,x_k\bm{]}$,
\begin{eqnarray}
\label{pt31.3}
\left( \widehat{f}_n(x) \,-\, f(x) \right)_+
& \leq & \left( \sup_{y\preceq x_k} \mbox{Av}_Y( \bm{(} y, \overline{x}_k \bm{]} ) \,-\, f(x) \right)_+
\nonumber \\
& \leq & \sup_{y\preceq x_k} \left| \mbox{Av}_\varepsilon( \bm{(} y, \overline{x}_k \bm{]} ) \right|
\,+\, \left( f(\overline{x}_k) \,-\, f(x) \right).
\end{eqnarray}
We apply Lemma~\ref{LA.1} to $\widetilde f(\widetilde x_1,\dots, \widetilde x_{d_2})= f( k_1,\ldots,k_{d_1},G_1^{-1}(\widetilde x_1),\ldots,G_{d_2}^{-1}(\widetilde x_{d_2}))$, $(\widetilde x_1,\dots, \widetilde x_{d_2})'\in [0,1]^{d_2},$ with $M=\widetilde M_n$. Since $\nu^d(\widetilde{B}_k)=\widetilde{h}_n^{d_2}=O(n^{-d_2/(d_2+2)})$ we obtain for the bias that
\begin{eqnarray}
\label{pt31.4}
\lefteqn{ \sum_{k\colon \, 1<k_{d_1+1},\ldots,k_d<\widetilde{M}_n}
\int_{\widetilde{B}_k} \left( f(\overline{x}_k) \,-\, f(x) \right) \, \nu(dx) } \nonumber \\
& \leq & \sum_{k\colon \, 1<k_{d_1+1},\ldots,k_d<\widetilde{M}_n} \left( f(\overline{x}_k) \,-\, f(\underline{x}_k) \right)
\; \nu(\widetilde{B}_k)
\,=\, O\left( \widetilde{M}_n^{d_2-1} \; \widetilde{h}_n^{d_2} \right)
\,=\, O\left( n^{-1/(d_2+2)} \right). \qquad
\end{eqnarray}

\noindent
We define again a dyadic scheme of nested hyperrectangles: For $j_1,\ldots,j_{d_2}\geq 0$,
\begin{displaymath}
\widetilde{B}_k^{(j_1,\ldots,j_{d_2})}
\,=\, \{(k_1,\ldots,k_{d_1})'\}\times
((k_{d_1+1}+1-2^{j_1})\widetilde{h}_n, (k_{d_1+1}+1)\widetilde{h}_n] \times
\cdots \times ((k_d+1-2^{j_{d_2}})\widetilde{h}_n, (k_d+1)\widetilde{h}_n]
\end{displaymath}
and
\begin{displaymath}
\widetilde{B}_{k,0}^{(j_1,\ldots,j_{d_2})}
\,=\, \{0,\ldots,k_1\}\times \cdots\times \{0,\ldots,k_{d_1}\} \times
((k_{d_1+1}+1-2^{j_1})\widetilde{h}_n, (k_{d_1+1}+1)\widetilde{h}_n] \times
\cdots \times ((k_d+1-2^{j_{d_2}})\widetilde{h}_n, (k_d+1)\widetilde{h}_n].
\end{displaymath}
Since $\{y\preceq x_k\}\subseteq \bigcup_{j_1,\ldots,j_{d_2}\geq 0} \widetilde{B}_{k,0}^{(j_1+1,\ldots,j_{d_2}+1)}\setminus \widetilde{B}_{k,0}^{(j_1,\ldots,j_{d_2})}$
and since $y\in \widetilde{B}_{k,0}^{(j_1+1,\ldots,j_{d_2}+1)}\setminus \widetilde{B}_{k,0}^{(j_1,\ldots,j_{d_2})}$ implies that
$\bm{(} y, \overline{x}_k \bm{]}\succeq \widetilde{B}_k^{(j_1,\ldots,j_{d_2})}$
we obtain, for all $x\in \widetilde B_k$,
\begin{displaymath}
\sup_{y\preceq x_k} \left| \mbox{Av}_\varepsilon( \bm{(} y, \overline{x}_k \bm{]} ) \right|
\,\leq\, \sum_{j_1,\ldots,j_{d_2}\geq 0}
\frac{ \sup\left\{ \left| \sum_{t=1}^n \varepsilon_t \1( I_t\in \bm{(} y,\overline{x}_k \bm{]} )
\right|\colon \;\; y\in \widetilde{B}_{k,0}^{(j_1+1,\ldots,j_{d_2}+1)}\setminus \widetilde{B}_{k,0}^{(j_1,\ldots,j_{d_2})} \right\} }
{ \# \left\{ t\leq n\colon \;\; I_t\in \widetilde{B}_k^{(j_1,\ldots,j_{d_2})} \right\} }.
\end{displaymath}
Recall that if the event~$\widetilde{A}_n$ occurs, then $\#\{t\leq n\colon \; I_t\in \widetilde{B}_k\}\geq (C_1/2)\, n^{2/(d_2+2)}$
for all~$k\in \widetilde{K}_n$, which implies that
\begin{displaymath}
\# \left\{t\leq n\colon \;\; I_t\in \widetilde{B}_k^{(j_1,\ldots,j_{d_2})} \right\} \,\geq\, (C_1/2) \; 2^{j_1+\cdots +j_{d_2}}
\; n^{2/(d_2+2)}.
\end{displaymath}
Furthermore, it follows from Lemma~\ref{LA.3} that for some $C<\infty$
\begin{equation}
\label{pt31.5}
E\left[ \sup\left\{ \left| \sum_{t=1}^n \varepsilon_t \1( I_t\in \bm{(} y,\overline{x}_k \bm{]} )
\right|\colon \;\; y\in \widetilde{B}_{k,0}^{(j_1+1,\ldots,j_{d_2}+1)} \right\} \right]
\,\leq\, C \; 2^{(j_1+\cdots +j_{d_2})/2} \; n^{1/(d_2+2)}.
\end{equation}
Therefore, we obtain that
\begin{eqnarray*}
\lefteqn{ \sup_{x\in \widetilde{B}_k} E\left[ \sup_{y\preceq x_k}
\left| \mbox{Av}_\varepsilon( \bm{(} y, \overline{x}_k \bm{]} ) \right| \;\; \1_{\widetilde{A}_n} \right] } \nonumber \\
& \leq & \sum_{j_1,\ldots,j_{d_2}\geq 0} E\left[
\frac{ \sup \left\{ \left| \sum_{t=1}^n \varepsilon_t \; \1(I_t\in \bm{(} y, \overline{x}_k \bm{]} ) \right|
\colon \quad y\in \widetilde{B}_{k,0}^{(j_1+1,\ldots,j_{d_2}+1)} \right\} }
{ \#\{t\leq n\colon\; I_t\in \widetilde{B}_k^{(j_1,\ldots,j_{d_2})}\} }
\;\; \1_{\widetilde{A}_n} \right] \nonumber \\
& = & O\left( n^{-1/(d_2+2)} \; \sum_{j_1,\ldots,j_{d_2}\geq 0} 2^{-(j_1+\cdots +j_{d_2})/2} \right)
\,=\, O\left( n^{-1/(d_2+2)} \right).
\end{eqnarray*}
This yields, in conjunction with (\ref{pt31.3}) and (\ref{pt31.4}), that (\ref{pt31.1}) holds.
The proof of (\ref{pt31.2}) is completely analogous and therefore it is omitted.
\end{proof}
\medskip

\subsection{Some auxiliary results}
\label{SS4.2}

{\lem
\label{LA.1}
Suppose that $f\colon\; [0,1]^d\rightarrow \R$ is isotonic
and let, for $M\in\N$ and $k=(k_1,\ldots,k_d)$, $x_k=(k_1/M,\ldots,k_d/M)$. Then
\begin{displaymath}
\sum_{k\colon\, 0<k_1,\ldots,k_d<M} \left( f(x_{k+\bm 1}) \,-\, f(x_{k-\bm 1}) \right)
\,\leq\, 2d\; M^{d-1}\; \left( f(1,\ldots,1) \,-\, f(0,\ldots,0) \right).
\end{displaymath}
}
\medskip

\begin{proof}[Proof of Lemma~\ref{LA.1}]
Let ${\mathcal I}_0=\{k\colon\; 0<k_1,\ldots,k_d<M \mbox{ and } k_j=1 \mbox{ for at least one }j \}$.
We estimate the sum by considering  the main and minor diagonals as follows:
\begin{eqnarray*}
\sum_{k\colon\; 0<k_1,\ldots,k_d<M} \left( f(x_{k+\mathbf{1}}) \,-\, f(x_{k-\mathbf{1}}) \right)
& = & \sum_{k\in {\mathcal I}_0} \sum_{i\geq 0} \left( f(x_{k+(i+1)\mathbf{1}}) \,-\, f(x_{k+(i-1)\mathbf{1}}) \right) \\
& \leq & \# {\mathcal I}_0 \;\; 2\left( \sup_x\{ f(x)\} \,-\, \inf_x \{ f(x)\} \right).
\end{eqnarray*}
The assertion of the lemma follows because   $\# {\mathcal I}_0\leq d M^{d-1}$.
\end{proof}
\medskip

{\lem
\label{LA.2}
Suppose that the assumptions of Theorem~\ref{T2.1} hold true.
Then, for arbitrary $\underline{z}\preceq \overline{z}$ with
$\bm[ \underline{z}, \overline{z} \bm] \subseteq [0,1]^d$ and some $\bar C<\infty$,
\begin{subequations}
\begin{equation}
\label{la2.1}
E\left[ \sup_{z\colon\; \underline{z}\preceq z\preceq\overline{z}}
\left| \frac{1}{\sqrt{n}} \sum_{t=1}^n \varepsilon_t \1\left( I_t\in \bm{(} z, \overline{z} \bm{]} \right) \right| \right]
\,\leq\, \bar C \; \sqrt{ P(I_1\in \bm{(} \underline{z}, \overline{z} \bm{]}) }
\end{equation}
and
\begin{equation}
\label{la2.2}
E\left[ \sup_{z\colon\; \underline{z}\preceq z\preceq\overline{z}}
\left| \frac{1}{\sqrt{n}} \sum_{t=1}^n \varepsilon_t \1\left( I_t\in \bm{[} \underline{z}, z \bm{)} \right) \right| \right]
\,\leq\,\bar C \; \sqrt{ P(I_1\in \bm{[} \underline{z}, \overline{z} \bm{)}) }.
\end{equation}
\end{subequations}
}
\medskip

\begin{proof}[Proof of Lemma~\ref{LA.2}]
We prove only (\ref{la2.1}) since the proof of (\ref{la2.2}) is completely analogous.
One of the main tools  which is used  is given by \citet[Thm. 1]{BW71}.
For this purpose, we adopt some notation from there.
A block~$B$ in~$\bm{[} \underline{z},\overline{z} \bm{]}$ is a subset of
$\bm{[} \underline{z},\overline{z} \bm{]}$ of the form
$\bm{(} u,v \bm{]} = (u_1,v_1]\times \cdots\times (u_d,v_d]$.
For $p\in\{1,\ldots,d\}$, the $p$th face of $B=\bm( u,v \bm]$ is
$(u_1,v_1]\times\cdots\times (u_{p-1},v_{p-1}]\times (u_{p+1},v_{p+1}]\times\cdots\times (u_d,v_d]$.
Disjoint blocks~$B$ and~$C$ are $p$-neighbors if they are abut and have the same $p$th face;
they are neighbors if they are $p$-neighbors for some~$p$.
(For example,
$(u_1,v_1]\times\cdots\times (u_{p-1},v_{p-1}]\times(\widetilde u,\widetilde v]\times (u_{p+1},v_{p+1}]\times\cdots\times (u_d,v_d]$
and
$(u_1,v_1]\times\cdots\times (u_{p-1},v_{p-1}]\times(\widetilde v,\widetilde w]\times (u_{p+1},v_{p+1}]\times\cdots\times (u_d,v_d]$
are $p$-neighbors if $0\leq \widetilde u<\widetilde v<\widetilde w\leq 1$.)
For each block $B=\bm{(} u,v \bm{]}$, let
\begin{displaymath}
X(B) \,=\, \frac{1}{\sqrt{n}} \sum_{t=1}^n \varepsilon_t \1( I_t\in B ) .
\end{displaymath}
\noindent
In what follows we show that condition~(2) in \citet[Thm.~1]{BW71} is fulfilled.
To this end, let~$B$ and~$C$ be arbitrary neighboring blocks in
$\bm{[} \underline{z},\overline{z} \bm{]}$. We will estimate the expected value
of the term
\begin{eqnarray*}
|X(B)|^2 |X(C)|^2 = \frac{1}{n^{2}} \sum_{t_1,t_2,t_3,t_4=1}^n \1(I_{t_1}\in B) \1(I_{t_2}\in B)
\1(I_{t_3}\in C) \1(I_{t_4}\in C) \varepsilon_{t_1} \varepsilon_{t_2} \varepsilon_{t_3} \varepsilon_{t_4}.
\end{eqnarray*}
Since~$B$ and~$C$ are disjoint sets it follows that
$$\1(I_{t_1}\in B) \1(I_{t_2}\in B) \1(I_{t_3}\in C) \1(I_{t_4}\in C)=0,$$
if $\{t_1,t_2\}\cap\{t_3,t_4\}\neq\emptyset$.
Therefore, and by independence of $(I_1',\varepsilon_1)',\ldots,(I_n',\varepsilon_n)'$,
\begin{eqnarray*}
\lefteqn{ E\left[ |X(B)|^2 |X(C)|^2 \right] } \\
& = & \frac{1}{n^2} \sum_{(t_1,\ldots,t_4)\colon \, \{t_1,t_2\}\cap\{t_3,t_4\}=\emptyset}
E[ \1(I_{t_1}\in B) \1(I_{t_2}\in B) \varepsilon_{t_1} \varepsilon_{t_2} ] \;
E[ \1(I_{t_3}\in C) \1(I_{t_4}\in C) \varepsilon_{t_3} \varepsilon_{t_4} ].
\end{eqnarray*}
Furthermore, again by independence of $(I_1',\varepsilon_1)',\ldots,(I_n',\varepsilon_n)'$,
and since
$$E[ \1(I_s\in B) \1(I_t\in B) \varepsilon_s \varepsilon_t ]
=E[ \1(I_s\in C) \1(I_t\in C) \varepsilon_s \varepsilon_t ]=0,$$
if $s\neq t$, we obtain that
\begin{eqnarray}
\label{pla2.1}
E\left[ |X(B)|^2 |X(C)|^2 \right]
& = & \frac{1}{n^2} \sum_{s\neq t} E\left[ \1(I_s\in B) \varepsilon_s^2 \right] \;
E\left[ \1(I_t\in C) \varepsilon_t^2 \right] \nonumber \\
& \leq & \overline{\sigma}_\varepsilon^4 \; P(I_1\in B) \; P(I_1\in C).
\end{eqnarray}
Let $m(B,C) = \min\{|X(B)|,|X(C)|\}$. From (\ref{pla2.1}) we obtain by Markov's inequality
that
\begin{displaymath}
P\left( m(B,C) \,\geq\, \nu \right)
\,\leq\, \frac{ E[ m(B,C)^4 ] }{ \nu^4 }
\,\leq\, \frac{ E[ |X(B)|^2 |X(C)|^2 ] }{ \nu^4 }
\,\leq\, \nu^{-4} \; \mu(B) \; \mu(C)
\end{displaymath}
for all $\nu>0$ and the measure $\mu(\cdot)=\overline{\sigma}_\varepsilon^2 P^{I_1}(\cdot)$.
Hence, condition~(2) of \citet[Thm. 1]{BW71} is fulfilled and it follows that
\begin{displaymath}
P\left( \sup_{z\colon\, \underline{z}\preceq z \preceq \overline{z}}
\left| \frac{1}{\sqrt{n}} \sum_{t=1}^n \varepsilon_t \1( I_t\in \bm (\underline{z}, z \bm] )
\right| \,\geq\, \nu \right)
\,\leq\, \widetilde C \; \nu^{-4} \; \mu( \bm{(} \underline{z}, \overline{z} \bm{]} )^2,
\end{displaymath}
for all $\nu>0$ and some $\widetilde C<\infty$.
This, however, implies that
\begin{eqnarray*}
\lefteqn{ E\left[ \sup_{z\colon\, \underline{z}\preceq z\preceq\overline{z}}
\left| \frac{1}{\sqrt{n}} \sum_{t=1}^n \varepsilon_t \1( I_t\in \bm{(} z, \overline{z} \bm{]} )
\right| \right] } \\
& = & \int_0^\infty P\left( \sup_{z\colon\, \underline{z}\preceq z\preceq\overline{z}}
\left| \frac{1}{\sqrt{n}} \sum_{t=1}^n \varepsilon_t \1( I_t\in \bm{(} z, \overline{z} \bm{]} )
\right| \,\geq\, \lambda \right) \, d\lambda \\
& \leq & \sqrt{ \mu( \bm{(} \underline{z}, \overline{z} \bm{]} ) }
\,+\, \int_{ \sqrt{\mu( \bm{(} \underline{z}, \overline{z} \bm{]} )}}^\infty
P\left( \sup_{z\colon\, \underline{z}\preceq z\preceq\overline{z}}
\left| \frac{1}{\sqrt{n}} \sum_{t=1}^n \varepsilon_t \1( I_t\in \bm{(} z, \overline{z} \bm{]} )
\right| \,\geq\, \lambda \right) \, d\lambda \\
& \leq & \sqrt{ \mu( \bm{(} \underline{z}, \overline{z} \bm{]} ) }
\,+\, 3\; \widetilde C \; \sqrt{ \mu( \bm{(} \underline{z}, \overline{z} \bm{]} ) },
\end{eqnarray*}
which proves the assertion of the lemma.
\end{proof}

{\lem
\label{LA.3}
Suppose that the assumptions of Theorem~\ref{T3.1} hold true.
Define $\rho_n=n^{-1}\sum_{t=1}^n P^{I_{n,t}}$.
Then, for arbitrary $\underline{z}\preceq \overline{z}$ with
$\bm[ \underline{z}, \overline{z} \bm] \subseteq \N_0^{d_1}\times \R^{d_2}$ and some $\bar C<\infty$,
\begin{subequations}
\begin{equation}
\label{la3.1}
E\left[ \sup_{z\colon\; \underline{z}\preceq z\preceq\overline{z}}
\left| \frac{1}{\sqrt{n}} \sum_{t=1}^n \varepsilon_t \1\left( I_t\in \bm{(} z, \overline{z} \bm{]} \right) \right| \right]
\,\leq\, \bar C \; \sqrt{ \rho_n\left( \bm{(} \underline{z}, \overline{z} \bm{]} \right) }
\end{equation}
and
\begin{equation}
\label{la3.2}
E\left[ \sup_{z\colon\; \underline{z}\preceq z\preceq\overline{z}}
\left| \frac{1}{\sqrt{n}} \sum_{t=1}^n \varepsilon_t \1\left( I_t\in \bm{[} \underline{z}, z \bm{)} \right) \right| \right]
\,\leq\, \bar C \; \sqrt{ \rho_n\left( \bm{[} \underline{z}, \overline{z} \bm{)} \right) }.
\end{equation}
\end{subequations}
}
\medskip

\begin{proof}[Proof of Lemma~\ref{LA.3}]
The proof is pretty much the same as that of Lemma~\ref{LA.2}.
Since we impose  condition (A3)(ii), we have only a bound for the conditional probability
$P(I_t\in C\mid I_1,\ldots,I_{t-d},\varepsilon_1,\ldots,\varepsilon_{t-d})$ but not for
$P(I_t\in C\mid I_1,\ldots,I_{t-1},\varepsilon_1,\ldots,\varepsilon_{t-1})$ at our disposal.
In view of this, we consider first the $d$-thinned partial sums
\begin{displaymath}
X_i(B) \,=\, \frac{1}{\sqrt{n}} \sum_{s\colon \; 1\leq sd+i\leq n}
\varepsilon_{sd+i} \; \1( I_{sd+i}\in B ),
\end{displaymath}
for $i=1,\ldots,d$, instead of  the full partial sums.
In analogy to (\ref{pla2.1}) in the proof of Lemma~\ref{LA.2} we show that,
for any neighboring blocks~$B$ and~$C$ in $\bm{[} \underline{z},\overline{z} \bm{]}$
and any $i\in\{1,\ldots,d\}$,
\begin{equation}
\label{pla3.1}
E\left[ |X_i(B)|^2 \; |X_i(C)|^2 \right]
\,\leq\, \widetilde C\; \rho_n(B) \; \rho_n(C),
\end{equation}
for some $\widetilde C<\infty$.

\noindent
As in the independent regressors case, we consider again, for arbitrary neighbored blocks~$B$ and~$C$
and arbitrary $t_1,t_2,t_3,t_4\in\{1,\ldots,n\}$,
the terms $E[\1(I_{t_1}\in B)\1(I_{t_2}\in B)\1(I_{t_3}\in C)\1(I_{t_4}\in C)
\varepsilon_{t_1}\varepsilon_{t_2}\varepsilon_{t_3}\varepsilon_{t_4}]$.
Since~$B$ and~$C$ are disjoint sets it follows as before that
$\1(I_{t_1}\in B)\1(I_{t_2}\in B)\1(I_{t_3}\in C)\1(I_{t_4}\in C)=0$ provided that
$\{t_1,t_2\}\cap\{t_3,t_4\}\neq\emptyset$. This  implies that the above
expectation is equal to 0.
Moreover, if the largest index appears only once, then the expectation also vanishes
since, by (A2)(i),
$E\left( \varepsilon_t \mid I_1,\ldots,I_t,\varepsilon_1,\ldots,\varepsilon_{t-1} \right)=0$.
Therefore, we have to examine in more detail  two cases:
$1\leq t_1,t_2<t_3=t_4\leq n$ and $1\leq t_3,t_4<t_1=t_2\leq n$.
Hence we obtain that
\begin{eqnarray*}
\lefteqn{ E\left[ |X_i(B)|^2 |X_i(C)|^2 \right] } \\
& = & \frac{1}{n} \sum_{t\colon\; d< td+i\leq n}
E\left[ \left( \frac{1}{\sqrt{n}} \sum_{s=0}^{t-1} \1(I_{sd+i}\in B) \varepsilon_{sd+i} \right)^2
\;\; \1(I_{td+i}\in C) \varepsilon_{td+i}^2 \right] \\
& & {} \,+\, \frac{1}{n} \sum_{t\colon\; d< td+i\leq n}
E\left[ \left( \frac{1}{\sqrt{n}} \sum_{s=0}^{t-1} \1(I_{sd+i}\in C) \varepsilon_{sd+i} \right)^2
\;\; \1(I_{td+i}\in B) \varepsilon_{td+i}^2 \right] \\
& \leq & \overline{\sigma}_\varepsilon^2 \;
\frac{1}{n} \sum_{t\colon\; d< td+i\leq n}  \;
E\left[ \left( \frac{1}{\sqrt{n}} \sum_{s=0}^{t-1} \1(I_{sd+i}\in B) \varepsilon_{sd+i} \right)^2\;P(I_{td+i}\in C\mid I_1,\ldots,I_{(t-1)d+i},\varepsilon_1,\ldots,\varepsilon_{(t-1)d+i}) \right] \\
& & {} \,+\, \overline{\sigma}_\varepsilon^2 \;
\frac{1}{n} \sum_{t\colon\; d< td+i\leq n}  \;
E\left[ \left( \frac{1}{\sqrt{n}} \sum_{s=0}^{t-1} \1(I_{sd+i}\in C) \varepsilon_{sd+i} \right)^2 \;P(I_{td+i}\in B\mid I_1,\ldots,I_{(t-1)d+i},\varepsilon_1,\ldots,\varepsilon_{(t-1)d+i}) \right] \\
& \leq & \widetilde C\; \rho_n(B) \; \rho_n(C),
\end{eqnarray*}
as required. Using (\ref{pla3.1}) to  estimate (\ref{pla2.1}), similar to  the proof of Lemma~\ref{LA.2},
we obtain in the same manner that
\begin{displaymath}
E\left[ \sup_{z\colon\; \underline{z}\preceq z\preceq\overline{z}}
\left| X_i\left( \bm{(} z, \overline{z} \bm{]} \right) \right| \right]
\,\leq\, \bar C \; \sqrt{ \rho_n\left( \bm{(} \underline{z}, \overline{z} \bm{]} \right) }.
\end{displaymath}
\noindent
Finally, summing up over $i=1,\ldots,d$ we obtain (\ref{la3.1}).
The proof of (\ref{la3.2}) is analogous and therefore it is omitted.
\end{proof}

\bigskip
\begin{ack}
This research was partly funded by the German Research Foundation DFG, project NE 606/2-2
and by the Volkswagen Foundation (Professorinnen für Niedersachsen des Nieders\"achsischen Vorab).
Part of this work was completed while K. Fokianos was visiting the Department of Statistics at TU Dortmund.
\end{ack}

\bibliographystyle{harvard}

\end{document}